\documentclass[12pt]{amsart}
\usepackage{amssymb}
\usepackage{amscd}
\usepackage{amsmath}

\newcommand{\nc}{\newcommand}
\nc{\cal}{\mathcal} %% Needed for LaTeX2e
\nc{\la}{\langle} \nc{\ra}{\rangle}

 \nc{\CA}{\cal A} \nc{\CBB}{\cal B}

 \nc{\CC}{\cal C}

 \nc{\CDD}{\cal D}
\nc{\CE}{\cal E}
\nc{\CF}{\cal F}
\nc{\CG}{\cal G}
\nc{\CH}{\cal H}
\nc{\CI}{\cal I}
\nc{\CJ}{\cal J}
\nc{\CK}{\cal K}
\nc{\CL}{\cal L}
\nc{\CM}{\cal M}
\nc{\CN}{\cal N}
\nc{\CO}{\cal O}
\nc{\CP}{\cal P}
\nc{\CQ}{\cal Q}
\nc{\CR}{\cal R}
\nc{\CS}{\cal S}
\nc{\CT}{\cal T}
\nc{\CU}{\cal U}
\nc{\CV}{\cal V}
\nc{\CW}{\cal W}
\nc{\CZ}{\cal Z}

\nc{\fa}{\mathfrak a}
\nc{\fg}{\mathfrak g}
\nc{\fk}{\mathfrak k}
\nc{\fh}{\mathfrak h}
\nc{\fm}{\mathfrak m}
\nc{\fn}{\mathfrak n}
\nc{\fA}{\mathfrak A}
\nc{\fC}{\mathfrak C}
\nc{\fI}{\mathfrak I}
\nc{\fL}{\mathfrak L}
\nc{\fS}{\mathfrak S}

\nc{\nen}{\newenvironment}
\nc{\ol}{\overline}
\nc{\ul}{\underline}
\nc{\lra}{\longrightarrow}
\nc{\lla}{\longleftarrow}
\nc{\Lra}{\Longrightarrow}
\nc{\Lla}{\Longleftarrow}
\nc{\Llra}{\Longleftrightarrow}
\nc{\hra}{\hookrightarrow}
\nc{\iso}{\overset{\sim}{\lra}}
\nc{\Hom}{\mathrm{Hom}}

\nc{\notebox}[1]{\noindent\fbox{\parbox{12.5cm}{\sf #1}}\\[8pt]}
\nc{\Thm}[1]{Theorem~\ref{#1}} \nc{\Prop}[1]{Proposition~\ref{#1}}
\nc{\Lem}[1]{Lemma~\ref{#1}} \nc{\Cor}[1]{Corollary~\ref{#1}}
\nc{\Conj}[1]{Conjecture~\ref{#1}} \nc{\Claim}[1]{Claim~\ref{#1}}
\nc{\Defn}[1]{ Definition~\ref{#1}} \nc{\Exa}[1]{Example~\ref{#1}}
\nc{\Rem}[1]{Remark~\ref{#1}} \nc{\Note}[1]{Note~\ref{#1}}

\nc{\marg}{\marginpar}

\nen{thm}[1]{\label{#1}{\bf Theorem.\ } \em}{}
\nen{prop}[1]{\label{#1}{\bf Proposition.\ } \em}{}
\nen{lem}[1]{\label{#1}{\bf Lemma.\ } \em}{}
\nen{klem}[1]{\label{#1}{\bf Key Lemma.\ } \em}{}
\nen{cor}[1]{\label{#1}{\bf Corollary.\ } \em}{}
\nen{conj}[1]{\label{#1}{\bf Conjecture.\ } \em}{}

\nen{claim}[1]{\label{#1}{\bf Claim.\ } \em}{}
\nen{defn}[1]{\label{#1}{\bf Definition.\ } }{}
\nen{exa}[1]{\label{#1}{\bf Example.\ } }{}

\nen{rem}[1]{\label{#1}{\em Remark.\ } }{}
\nen{exer}[1]{\label{#1}{\em Exercise.\ } }{}
\nen{pf}{\begin{proof}}{\end{proof}}

\nc{\br}{\mathbb R} \nc{\bz}{\mathbb Z} \nc{\bc}{\mathbb C}
\nc{\bn}{\mathbb N} \nc{\geg}{\mathfrak g} \nc{\gh}{{\goth h}}
\nc{\gan}{{\goth a}} \nc{\G}{\Gamma} \nc{\g}{\gamma}
\nc{\sm}{\setminus} \nc{\sub}{\subset} \nc{\lm}{\lambda}
\nc{\eps}{\varepsilon} \nc{\nty}{\infty} \nc{\al}{\alpha}
\nc{\bt}{\beta} \nc{\om}{\omega} \nc{\dl}{\delta} \nc{\Dl}{\Delta}
\nc{\Om}{\Omega} \nc{\s}{\sigma} \nc{\ro}{\rho} \nc{\te}{\theta}
\nc{\SLR}{SL_2(\br)} \nc{\GLR}{GL_2(\br)} \nc{\PGLR}{PGL_2(\br)}
\nc{\PSLR}{PSL_2(\br)} \nc{\SLC}{SL(2,\bc)} \nc{\uH}{\mathbb H}
\nc{\fD}{\mathfrak D} \nc{\fE}{\mathfrak E} \nc{\haf}{\frac{1}{2}}
\nc{\qtr}{\frac{1}{4}} \nc{\8}{\infty} \nc{\7}{{-\infty}}
\nc{\inv}{^{-1}}
\begin{document}
%\topmatter
\title[automorphic functions]
{Estimates of automorphic functions}

%June 1, 2003
%to appear in special issue of Moscow Math. J. in honor of Pierre
%Cartier

\author{Joseph Bernstein and Andre Reznikov}
\address{Tel Aviv University, Ramat Aviv, Israel}
\email{bernstei@math.tau.ac.il}
\address{Bar Ilan University, Ramat-Gan, Israel}
\email{reznikov@math.biu.ac.il}

\begin{abstract}
We present a new method of estimating trilinear period for
automorphic representations of $\SLR$. The method is based on the
uniqueness principle in representation theory.
  We show how to separate the exponentially decaying factor in the
triple period from the essential automorphic factor which behaves
polynomially. We also describe a general method which gives an
estimate on the average of the automorphic factor and thus prove a
convexity bound for the triple period.
\end{abstract}
\maketitle

\section{Introduction}
\label{intro}
\subsection{Maass forms}
Let $Y$ be a compact Riemann surface  with a Riemannian metric of
constant curvature $-1$ and the associated volume element $dv$.
The corresponding  Laplace-Beltrami operator  is non-negative and
has purely discrete spectrum on the space $L^2(Y,dv)$ of functions
on $Y$. We will denote by $0=\mu_0< \mu_1 \leq \mu_2 \leq ...$ its
eigenvalues  and by $\phi_i=\phi_{\mu_i}$ the corresponding
eigenfunctions (normalized to have $L^2$ norm one). In the theory
of automorphic forms the functions $\phi_{\mu_i}$ are called
automorphic functions or {\bf Maass forms} (after H. Maass,
\cite{M}).

The study of Maass forms plays an important role in analytic
number theory.

 We are interested in their analytic properties and will present
a new method of bounding some important quantities arising from
$\phi_{i}$.

\subsection{Triple products} For any three Maass forms
$\phi_i,\ \phi_j, \ \phi_k$ we define the following
triple product or triple period:
\begin{eqnarray}\label{cijk}
c_{ijk}=\int_Y\phi_i\phi_j\phi_kdv \ .
\end{eqnarray}

We would like to bound the coefficient $c_{ijk}$ as a function of
eigenvalues  $\mu_i,\ \mu_j,\ \mu_k$. In particular, we would like
to find bounds for these coefficients  when one or more of these
indices tend to infinity.

\subsection{Motivation} First of all we would like to explain why
this problem is interesting. The explanation goes back to
pioneering works of Rankin and Selberg (see \cite{Ra}, \cite{Se}).
They discovered that in special cases triple products as above
give rise to automorphic $L$-functions. That allowed them to
obtain analytic continuation and effective bounds for these
$L$-functions and as an application to obtain bounds on Fourier
coefficients of cusp forms towards Ramanujan conjecture.

 Since then the Rankin-Selberg method has had many generalizations.
Recently, for $Y$ arising from the full modular group $SL_2(\bz)$
and for cuspidal functions $\phi$,
 Watson (see \cite{Wa}) proved the following beautiful formula:

\begin{eqnarray}\label{wats}
\left|\int_Y\phi_i\phi_j\phi_kdv\right|^2=
G(\lm_i,\lm_j,\lm_k)\frac{L(1/2,\phi_i\otimes\phi_j\otimes\phi_k)}
{ L(1,\phi_i,Ad)L(1,\phi_j,Ad)L(1,\phi_k,Ad)}\ .
\end{eqnarray}

Here $\lm_t$ is a natural parameter of an eigenfunction $\phi_t$
related to the eigenvalue by $\mu_t=\frac{1-\lm_t^2}{4}$. The
unctions $L(s,\phi_i\otimes\phi_j\otimes\phi_k)$ and
$L(s,\phi_t,Ad)$ are appropriate automorphic $L$-functions
associated to $\phi_i$, and the function $G(\lm_i,\lm_j,\lm_k)$ is
an explicit rational expression in the ordinary $\G$-functions.
The relation (\ref{wats}) can be viewed as a far reaching
generalization of the original Rankin-Selberg formula. It was
motivated by a work \cite{HK} by Harris and Kudla  on a conjecture
of Jacquet.

\subsection{Results}
In this paper we will consider the following problem. We fix two
Maass forms $\phi = \phi_\tau,\ \phi'= \phi_{\tau'}$ as above and
consider coefficients defined by the triple period as above:
\begin{eqnarray}\label{ci}
c_i=\int_Y\phi\phi'\phi_idv
\end{eqnarray}
as  $\{\phi_i = \phi_{\lm_i}\}$ run over the basis of Maass forms.

Thus we see from (\ref{wats}) that the estimates of coefficients
$c_i$ are equivalent to the estimates of the corresponding
$L$-functions. One would like to have a general method of
estimating the coefficients $c_i$ and similar quantities. This
problem was raised by Selberg in his celebrated paper \cite{Se}.

Let us  understand what kind of bounds on the left hand side of
(\ref{wats}) one would like to have in order to estimate
effectively $L$-functions involved in the right hand side of
(\ref{wats}) (or at least the ratio of $L$-functions).

We note first that one expects that $c_i$ have exponential decay
in $|\lm_i|$ as $i$ goes to $\8$. Namely, general experience from
the analytic theory of automorphic $L$-functions tells us that
$L$-functions have at most polynomial growth when  $|\lm_i|\to\8$.
Hence, analyzing the function $G(\lm)$,
 one would expect from (\ref{wats}) and
the Stirling formula for the asymptotics of $\G$-function that the
normalized coefficients
\begin{eqnarray}\label{ai}
b_i=|c_i|^2\exp(\frac{\pi}{2}|\lm_i|)
\end{eqnarray}
have at most polynomial growth in $|\lm_i|$, and hence $c_i$ decay
exponentially. However, it is difficult to see from the definition
of the coefficients $c_i$ that they have exponential decay and it
is not clear what should be the rate of this decay.

 The fact that an exponential decay with the exponent $\frac{\pi}{2}$ holds
 for a general Riemann surface was first shown by Good
 and Sarnak (see \cite {Go} and \cite{Sa}). Both proofs used
ingenious analytic continuation of automorphic functions in the
{\it variable} parameter.

In this paper we will explain how to naturally separate  the
exponential decay from a polynomial growth in coefficients $c_i$
using representation theory. We also prove the following

\begin{thm}{THM}
There exists an effectively computable
 constant $A$ such that the following bound holds for arbitrary $T > 0$

\begin{eqnarray}\label{mean-value}
\sum_{T\leq |\lm_i|\leq 2T}b_i\leq A\ .
\end{eqnarray}
\end{thm}

\subsection{A conjecture} \label{conj} The estimate in the
theorem is tight but if we try to use it to get a bound for an
individual term $b_i$  we get only an inequality
\begin{eqnarray}\label{basic}
b_i\leq A \ .
\end{eqnarray}

    According to Weyl's law there are approximately $c T^2$
eigenvalues $\mu_i$ with $\lm_i$ between $T$ and $2T$, so the
individual bound for the coefficient $b_i$ is definitely not
tight. We would like to make the following conjecture concerning
the size of coefficients $b_i$:

\begin{conj}{CONJ1}
For any $\eps>0$ there exists a constant $C_\eps>0$   such that
$$b_i\leq C_\eps|\lm_i|^{-2+\eps}\ ,$$
as $|\lm_i|\to\8$.
\end{conj}

For $Y$ arising from congruence subgroups this conjecture is
consistent with the Lindel\"{o}f conjecture for appropriate
automorphic $L$-functions (see \cite{BR1}, \cite{Sa2} and
\cite{Wa} for more details). We note that the bound in the Theorem
above corresponds to the so--called {\bf convexity bound}.

\subsection {The method} The first proof of the (slightly weaker)
version of Theorem \ref{THM} appeared in \cite{BR1}. It was based
on the analytic continuation of representations from a real group
to a complex group (generalizing methods of \cite{Sa}). The method
based on the analytic continuation was extended in \cite{KS} to
the case of
 higher rank groups. While it gives bounds which are tight for
general representations, it was not able, so far, to cover cases
relevant to $L$-functions.

The proof we present here is based on the uniqueness of triple
product in representation theory. It has an advantage that it
could be generalized to higher rank groups and gives bounds which
are consistent with the theory of $L$-functions.  The present
method also could be applied  to $p$-adic groups (unlike methods
of \cite{BR1}).

%%%%%%%%%%%%%%%%%%%%
We describe now the general ideas behind our new proof. It is
based on ideas from representation theory. Namely, we use the fact
that every automorphic form $\phi$ generates an automorphic
representation of the group $G = \PGLR$; this means that starting
from $\phi$ we produce a smooth irreducible representation of the
group $G$ in a space $V$ and its realization $\nu : V \to
C^{\8}(X)$ in the space of smooth functions on the automorphic
space $X = \G \backslash G$.

 The triple product  $c_i=\int_Y\phi\phi'\phi_idv$ extends to a
$G$-equivariant trilinear form on the corresponding automorphic
representations  $l^{aut}:V\otimes V'\otimes V_i\to\bc$, where $V
= V_\tau, V' = V_{\tau'}, V_i = V_{\lambda_i}$ .

Then we use a general result from representation theory that such
$G$-equivariant trilinear form is unique up to a scalar. This
implies that the automorphic form $l^{aut}$ is proportional to an
explicit "model" form $l^{mod}$ which we describe using explicit
realizations of representations of the group $G$; it is important
that this last form carries no arithmetic information.

Thus we can write $l^{aut} = a_i \cdot l^{mod}$ for some constant
$a_i$ and hence $c_i= l^{aut}(e_\tau\otimes e_{\tau'} \otimes
e_{\lambda_i}) = a_i \cdot l^{mod}(e_\tau \otimes e_{\tau'}
\otimes e_{\lambda_i})$, where $e_\tau,\ e_{\tau'},\
e_{\lambda_i}$ are K-invariant unit vectors in the automorphic
representations $V, V',\ V_i$ corresponding to the automorphic
forms $\phi$, $\phi'$ and $\phi_i$.

 It turns out that the proportionality coefficient
 $ a_i$ in the last formula carries an
 important "automorphic" information while the second factor
 carries no arithmetic information and can be computed
 in terms of $\G$-functions using explicit realizations of
 representations $V_\tau$, $V_{\tau'}$ and $V_{\lambda_i}$. This second factor is
 responsible for
 the exponential decay, while the first factor $a_i$ has a
 polynomial behavior in parameter $\lm_i$.

   In order to bound the quantities $a_i$, we use the fact that they
   appear as coefficients in the spectral
   decomposition of the diagonal Hermitian form $H_{\Dl}$
   on the space $E = V_\tau \otimes V_{\tau'}$ (see \ref{diagonal}, \ref{first}).
   This gives an inequality $\sum |a_i|^2 H_i \leq H_\Dl$
   where $H_i$ is an Hermitian
form on $E$ induced by the model trilinear form $l^{mod}: V
\otimes V' \otimes V_i \to \bc$ as above.

    Using the geometric properties of
   the diagonal form and simple explicit estimates of forms $H_i$
   we establish the convexity bound for the coefficients $a_i$.

  It is known that the uniqueness principle plays a central role in
  the theory of automorphic functions (see \cite{PS}).
  The impact that the uniqueness has on
the analytic behavior of automorphic functions is yet another
manifestation of this principle.

%%%%%%%%%%%%%%%%%%%%%%%%%%%%%%%%%%%%%%%%%%%%%%%%%%%%%%%%%

{\bf Acknowledgments.}  We would like to thank Peter Sarnak for
many fruitful discussions. We also would like to thank Yuri
Neretin for discussions on trilinear functionals and the referee
for helpful comments.

Research was partially supported by EC TMR network "Algebraic Lie
Representations", grant no.  ERB FMRX-CT97-0100, by BSF grant,
Minerva Foundation and by the Excellency Center ``Group Theoretic
Methods in the Study of Algebraic Varieties'' of the  Israel
Science Foundation, the Emmy Noether Institute for Mathematics
(the Center of Minerva Foundation of Germany).

The paper was written during our stay at IHES and is an extended
version of the talk given by the first author at Cartier
Colloquium in June 2002.

%%%%%%%%%%%%%%%%%%%%%%%%%%%%%%%%%%%%%%%%%%%%%%%%%%%%%%%%%%
\section{Representation theoretic setting}
\label{reps} We recall the standard connection of the above
setting with  representation theory (see \cite{GGPS}).

\subsection{Automorphic functions and automorphic representations}
\label{ureps}

Let us describe the geometric construction which allows one to
pass from analysis on a Riemann surface to representation theory.

Let $\uH$ be the upper half plane with the hyperbolic metric of
constant curvature $-1$. The group $\SLR$ acts on $\uH$ by
fractional linear transformations. This action  allows to identify
the group $\PSLR$ with the group of all orientation preserving
motions of $\uH$. For reasons explained bellow we would like to
work with the group $G$ of all motions of $\uH$; this group  is
isomorphic to $\PGLR$. Hence throughout the paper we denote
$G=\PGLR$.

Let us fix a discrete co-compact subgroup $\G \subset G$ and set
$Y=\G \sm \uH$. We consider the Laplace operator on the Riemann
surface $Y$ and denote by $\mu_i$ its eigenvalues and by $\phi_i$
the corresponding normalized eigenfunctions.

  The case when $\G$ acts freely on $\uH$ precisely corresponds to the
case discussed in the introduction (this follows from the
uniformization theorem for the Riemann surface $Y$). Our results
hold for general co-compact subgroup $\G$ (and in fact, with
slight modifications, for any lattice $\G \subset G$).

    We will identify the upper half plane $\uH$ with $G / K$, where
$K = PO(2)$ is a maximal compact subgroup of $G$ ( this follows
from the fact that $G$ acts transitively on $\uH$ and the
stabilizer in $G$ of the point
 $z_0 = i \in \uH$ coincides with $K$).

We denote by $X$ the compact quotient $\G\sm G$ (we call it the
automorphic space). In the case when  $\G$ acts freely  on $\uH$
one can identify the space $X$ with the bundle of unit tangent
vectors to the  Riemann surface $Y = \G \sm \uH$.

 The group $G$ acts on $X$ (from the right) and hence on the space of
 functions on $X$.
 We fix the unique $G$-invariant measure $\mu_X$ on $X$ of total
mass one. Let $L^2(X)=L^2(X,d\mu_X)$ be the space of square
integrable functions and $(\Pi_X, G, L^2(X))$ the corresponding
unitary representation. We will denote by $P_X$ the Hermitian form
on $L^2(X)$ given by the scalar product. We denote by $||\ ||_{X}$
or simply $||\ ||$ the corresponding norm and by $\langle f,g
\rangle_X$ the corresponding scalar product.

The identification  $Y=\G\sm \uH\simeq X/K$ induces the embedding
$L^2(Y)\sub L^2(X)$.
 We will always identify the space $L^2(Y)$ with the subspace of
$K$-invariant functions in  $L^2(X)$.

Let  $\phi$ be a normalized eigenfunction of the Laplace-Beltrami
 operator on $Y$. Consider a closed $G$-invariant subspace
$L_\phi\sub L^2(X)$ generated by $\phi$ under the  action of $G$.
It is well-known that $(\pi,L)=(\pi_\phi, L_\phi)$ is an
irreducible unitary representation of $G$ (see \cite{GGPS}).

    Usually it is more convenient to work with the space $V =
    L^\8$ of smooth vectors in $L$. The unitary Hermitian form $P_X$
    on $V$ is $G$-invariant.

    A smooth representation $(\pi, G, V)$ equipped with
    a positive $G$-invariant Hermitian form $P$ we will call
   a {\bf smooth pre-unitary representation}; this simply means that
    $V$ is the space of smooth vectors in the unitary
    representation obtained from $V$ by completion with respect to
    $P$.

   Thus starting with an automorphic function $\phi$ we constructed
   an irreducible smooth pre-unitary representation $(\pi, V)$.
   In fact we constructed this space together with a canonical
   morphism $\nu : V  \to C^\8 (X)$ since $C^\8(X)$ is the smooth
   part of $L^2(X)$.

   \defn {enhanced}A smooth pre-unitary representation
   $(\pi, G, V)$ equipped with a $G$-morphism $\nu: V \to C^\8(X)$
   we will call an {\bf $X$-enhanced representation}.

   In this note we will assume that the morphism $\nu$ is
   normalized,
    i.e. it carries the standard $L^2$ Hermitian
   form $P_X$ on $C^\8(X)$ into Hermitian form $P$ on $V$.

   Thus starting with an automorphic function $\phi$ we
   constructed

   (i) An $X$-enhanced irreducible pre-unitary representation
   $(\pi, V, \nu)$,

   (ii) A $K$-invariant unit vector $e_V \in V$
   (this vector is just our function $\phi$).

   Conversely, suppose we are given an irreducible smooth
   pre-unitary
    $X$-enhanced representation $(\pi,V, \nu)$ of
the group $G$ and a $K$-fixed unit vector $e_V \in V$. Then the
function $\phi = \nu(e_V) \in C^\8(X)$ is $K$-invariant and hence
can be considered as a function on $Y$. The fact that the
representation $(\pi, V)$ is irreducible implies that $\phi$ is an
automorphic function.

 Thus we have established a natural correspondence between
 Maass forms $\phi$ and tuples $(\pi, V, \nu,
 e_V)$,   where $(\pi, V,\nu)$ is an $X$-enhanced irreducible
 smooth pre-unitary representation and $e_V \in V$ is a unit
 $K$-invariant vector.

\subsection{Decomposition of the representation $(\Pi_X, G, L^2(X))$}

   It is well known that in case when $X$ is compact the
   representation $(\Pi_X, G, L^2(X))$ decomposes into a direct
   (infinite) sum
\begin{equation}\label{spec-L2X}
L^2(X)=\oplus_j (\pi_j, L_j)
\end{equation}
of irreducible unitary representations of $G$ (all representations
appear with finite multiplicities (see \cite{GGPS})). Let $(\pi,
L)$ be one of these irreducible "automorphic" representations and
$V = L^\8$ its smooth part. By definition $V$ is given with a
$G$-equivariant isometric morphism $\nu: V \to C^\8(X)$, i.e. $V$
is an $X$-enhanced representation.

   If $V$ has a $K$-invariant vector it corresponds to a Maass form.
There are  other spaces in this decomposition which  correspond to
discrete series representations. Since they are not related to
Maass forms we will not study them in more detail.

\subsection{Representations of $\PGLR$}\label{irrrep}
 All irreducible
unitary representations of $G$ are classified. For simplicity we
consider those with a nonzero $K$-fixed vector (so called
representations of class one) since only these representations
arise from Maass forms. These are the representations of the
principal and the complementary series and the trivial
representation.

   We will use the following standard explicit model for irreducible
smooth   representations of  $G$.

    For every complex number $\lm$ consider
 the space $V_\lm$ of smooth even homogeneous functions on
$\br^2\sm0$ of homogeneous degree $\lm-1 \ $ (which means that
$f(ax,ay)=|a|^{\lm-1}f(x,y)$ for all $a\in\br \setminus 0$).
 The representation $(\pi_\lm, V_\lm)$ is induced by
the action of the group $\GLR$ given by $\pi_\lm (g)f(x,y)=
f(g\inv (x,y))|\det g|^{(\lm-1)/2}$. This action is trivial on the
center of $\GLR$ and hence defines a representation of $G$. The
representation $(\pi_\lm, V_\lm)$ is called {\bf representation of
the generalized principal series}.

When $\lm=it$ is  purely imaginary  the representation $(\pi_\lm
,V_\lm)$ is pre-unitary; the $G$-invariant scalar product in
$V_\lm$ is given by $\langle f,g \rangle_{\pi_\lm}=\frac{1}{2\pi}
\int_{S^1} f\bar g d\te$. These representations are   called
representations of {\bf the principal series}.

 When $\lm\in (-1,1)$ the representation $(\pi_\lm ,V_\lm)$ is called
 a representation of the complementary series. These representations
 are also pre-unitary, but the formula for the scalar product is
 more complicated (see \cite{GGPS}).

 All these representations have $K$-invariant vectors.
 We fix a $K$-invariant unit vector $e_{\lm} \in V_\lm$ to be
  a function which is one on the unit circle in $\br^2$.

  Representations of the principal and the complimentary series exhaust
  all nontrivial irreducible pre-unitary representations of $G$
  of class one.

 In what follows we will do necessary computations for
 representation of the principal series. Computations for the
 complementary series are a little more involved but essentially
 the same (compare with \cite{BR1}, section 5.5, where similar
 computations are described in detail).

Suppose we are given a class one $X$-enhanced representation
 $\nu: V_{\lm} \to C^\8(X)$; we assume $\nu$ to be an isometric embedding.
 Such $\nu$ gives rise to an
eigenfunction of the Laplacian on the Riemann surface $Y = X/K$ as
before. Namely, if $e_{\lm} \in V_\lm$ is a unit  $K$-fixed vector
then the function $\phi = \nu(e_\lm)$ is a normalized
eigenfunction of the Laplacian on the space $Y = X/K$ with the
eigenvalue $\mu=\frac{1-\lm^2}{4}$. This explains why $\lm$ is a
natural parameter to describe Maass forms.

\subsection{Triple products}\label{3prod-def}
We introduce now our main tool.
\subsubsection{Automorphic triple products}\label{aut3prod}
Suppose we are given three $X$-enhanced representations of $G$
\begin{eqnarray*}
\nu_j:V_j\to C^\8(X),\ \
 j=1,2,3\ .
\end{eqnarray*}

We define the  $G$-invariant trilinear form
$l^{aut}_{\pi_1,\pi_2,\pi_3}:V_1\otimes V_2\otimes V_3\to\bc$ ,
 by formula

\begin{eqnarray}\label{aut3prod-def}
l^{aut}_{\pi_1,\pi_2,\pi_3}(v_1\otimes v_2\otimes v_3)= \int_{
X}\phi_{v_1}(x)\phi_{v_2}(x)\phi_{v_3}(x)d\mu_X\ ,
\end{eqnarray}
where $\phi_{v_j}=\nu_j(v_j)\in C^\8(X)$ for $v_j\in V_j$.

In particular, the triple periods $c_i$  in (\ref{ci}) can be
expressed in terms of this form as $c_i=
l^{aut}_{\pi,\pi',\pi_i}(e_\tau\otimes e_{\tau'}\otimes
e_{\lm_i})$, where $e_\lm \in V_\lm$ is the  $K$-fixed unit
vector.

\subsubsection{Uniqueness of triple products} The central fact about
invariant trilinear functionals is the following uniqueness result:

\begin{thm}{ubi} Let $(\pi_j,V_j),\ \ j=1,2,3\ ,$
be three irreducible smooth admissible
  representations of  $G$. Then
$\dim\Hom_G(V_1\otimes V_2\otimes V_3,\bc)\leq 1$.
\end{thm}

\begin{rem}{rem-uni}
The  uniqueness statement was proven by Oksak in \cite{O} for the
group $\SLC$ and the proof could be adopted for $\PGLR$ as well
(see also \cite{Mo} and \cite{Lo}). For the $p$-adic $GL(2)$ more
refined results were obtained by Prasad (see \cite{Pr}). He also
proved the uniqueness when at least one representation is a
discrete series representation of $\GLR$.

 There is no uniqueness of trilinear functionals for representations
of $\SLR$ (the space is two-dimensional). This is the reason why we
prefer to work with $\PGLR$.

For $\SLR$ one has the following uniqueness statement instead. Let
$(\pi,V)$ and $(\sigma, W)$ be two  irreducible smooth pre-unitary
representations of  $\SLR$ of class one. Then the space of
$\SLR$-invariant trilinear functionals on $V\otimes V\otimes W$
which are symmetric in the first two variables is one-dimensional.
This is the correct uniqueness result needed if one wants to work
with $\SLR$ ; this was implicitly done in \cite{Re2}, where the
second author missed the absence of the uniqueness for $\SLR$. We
take an opportunity to correct this gap.

We note however, that the absence of uniqueness does not pose any
problem for the method we present. All what is really needed for
our method is the fact that the  space of invariant functionals is
finite dimensional .
\end{rem}

\section{Triple products: exponential decay}

We now explain our method how to bound coefficients $c_i$. It is
based on the uniqueness of trilinear functionals.

\subsection{Model triple products}\label{mod3prod-def}

   Let $(\pi, V)$ and $(\pi', V')$ be
   automorphic representations corresponding to Maass forms $\phi$
   and $\phi'$. Any Maass form $\phi_i$ gives us an automorphic
   representation $(\pi_i,V_{\lm_i})$ and hence defines a trilinear
   functional
   $$l^{aut}_{\pi,\pi',\pi_i}: V \otimes V' \otimes V_{\lm_i} \to
   \bc\ .$$

In \ref{modfunc} we use an explicit  model for representations
$\pi_1,\ \pi_2,\ \pi_3$ to construct a model invariant trilinear
functional which is given  by an explicit formula. We call it the
{\bf model triple product} and denote  by
$l^{mod}_{\pi_1,\pi_2,\pi_3}$.

By the uniqueness principle for representations $\pi,\pi',\pi_i$
there exists a constant $a_i=a_{\pi,\pi',\pi_i}$ such that:
\begin{eqnarray}\label{coef-a-def}
l^{aut}_{\pi,\pi',\pi_i}=a_i \cdot l^{mod}_{\pi,\pi',\pi_i}\ .
\end{eqnarray}

\subsection{Exponential decay}\label{expdecay}
   This gives a formula  for the triple products $c_i$
\begin{eqnarray}\label{aut-mod-vect}
c_i=l^{aut}_{\lm_i}(e_\tau\otimes e_{\tau'}\otimes e_{\lm_i})=a_i
\cdot l^{mod}_{\lm_i}(e_\tau\otimes e_{\tau'}\otimes e_{\lm_i}) \
.
\end{eqnarray}
Here we denoted $l^{aut}_{\lm_i}=l^{aut}_{\pi,\pi',\pi_i}$,
$l^{mod}_{\lm_i}=l^{mod}_{\pi,\pi',\pi_i}$ and $e_\lm$ is the unit
$K$-fixed vector in the representation $V_\lm$.

The model triple product $l^{mod}_{\lm_i}(e_\tau\otimes
e_{\tau'}\otimes e_{\lm_i})$  constructed in \ref{modfunc} is
given by an explicit integral. In Appendix \ref{A} we evaluate
this integral by a direct computation in the model. It turns out
that it has  an exponential decay in $|\lm|$ which explains the
exponential decay of coefficients $c_i$.
 Namely, we prove the following

\begin{prop}{propklambda} Set $k_\lm := |l^{mod}_{\lm}(e_\tau\otimes
e_{\tau'}\otimes e_{\lm})|^2$. Then there exists a constant $c>0$
such that
$$k_{\lm} = c\exp(-\frac{\pi}{2}|\lm|)\cdot |\lm|^{-2} (1 +
O(|\lm|^{-1})) $$
 as $|\lm|\to\8$ and $\lm\in i\br$.
\end{prop}

\section{Triple products: polynomial bounds}\label{pbound}

We explain now how to obtain bounds on the coefficients $a_i$
(note that these coefficients  encode deep arithmetic information
 - values of $L$-functions).

Our method is based on the fact that these coefficients appear in
the spectral decomposition of some geometrically defined Hermitian
form  on the space $E$ which is essentially the tensor product of
spaces $V$ and $V'$.

  More precisely, denote by $L$ and $L'$ the Hilbert completions
of spaces $V$ and $V'$, consider the unitary representation $(\Pi,
G \times G, L \otimes L')$ of the group $G \times G$ and denote by
$E$ its smooth part; so $E$ is a smooth completion of $V \otimes
V'$.

   Denote by $\CH(E)$ the  (real) vector space of continuous Hermitian
   forms on $E$
and by $\CH^+(E)$ the cone of nonnegative Hermitian forms.

   We will describe several classes of Hermitian forms on $E$;
   some of them have spectral description, others are described
   geometrically.

\subsection{Hermitian forms corresponding to trilinear
functionals} \label{diag-aut}

Let $W$ be a smooth pre-unitary admissible representation of $G$.
Any $G$-invariant functional $l: V\otimes V'\otimes W \to \bc$
defines a $G$-intertwining morphism $T^l: V \otimes V' \to W^*$
which extends to a $G$-morphism
\begin{eqnarray}\label{T-W-def}
T^l: E \to  \bar W \ .
\end{eqnarray}

where we identify the complex conjugate space $\bar W$   with the
smooth part of the space $W^*$.

   The standard Hermitian form (scalar product) $P_W$ on the space
$W$ induces the Hermitian form $\bar P$ on $\bar W$.
 Using the operator $T^l$ we
define the Hermitian form $H^l$ on the space $E$ by $H^l =
(T^l)^{\ast} (\bar P)$, i.e. $H^l(u)= \bar P(T^l(u))$  for $u\in
E$.

We note that if the representation of $G$ in the space $W$ is
irreducible then starting with the Hermitian form $H^l$ we can
reconstruct the space $W$,
 the functional $l$ and the morphism $T^l$ uniquely up to an
isomorphism.

Let us introduce a special notation for the particular case we are
interested in. For any number $\lm \in i \br$ consider the
representation of the principal series
  $W = V_\lm$, choose the model trilinear functional
  $l^{mod}: V \otimes V' \otimes V_\lm \to \bc$
  described in \ref{modfunc} and denote the
 corresponding Hermitian form on $E$ by $H^{mod}_\lm$.

\subsection{Diagonal form $H_\Dl$} \label{diagonal}

Consider the space $C^\8(X\times X)$. The diagonal $\Dl:X \to
X\times X$ gives rise to the restriction morphism $r_\Dl :
C^\8(X\times X)\to C^\8(X)$.
 We define a nonnegative Hermitian form $H_\Dl$ on $C^\8(X\times X)$ by
$H_\Dl = (r_\Dl)^{\ast}(P_X)$, i.e.

$H_\Dl(u)= P_X(r_\Dl(u)) = \int_{X}|r_\Dl(u)|^2d\mu_X $
 for $u\in C^\8(X\times X)$.

 We call $H_\Dl$ the diagonal form.

More generally, if $L$ is a closed subspace of $ L^2(X)$ and
$pr_L: L^2(X) \to L$ the orthogonal projection onto $L$ we can
define a Hermitian form   $P_L$ on $C^\8(X)$  by $P_L =
(pr_L)^{\ast} (P_X)$ and consider the induced Hermitian form $H_L=
(r_\Dl)^{\ast}(P_L)$ on $C^\8(X \times X)$.

   Clearly the correspondence $L \mapsto H_L$ is additive
   (which means that $H_{L + L'} = H_L + H_{L'}$ if $L$ and $L'$ are orthogonal)
   and monotone.

\subsection{First basic inequality }\label {first}
 Let us realize the space $E = V\otimes V'$ as a $G\times
G$-invariant
 subspace of $C^\8(X\times X)$. We consider the restrictions of
 the Hermitian forms $H_\Dl , H_L$ discussed above to the space
 $E$ and will denote them by the same symbols.

 \begin{claim}{claimfirst} Let $\phi_{\lm_i}$ be a Maass form. Consider the
 $G$-invariant subspace $L_i \subset L^2(X)$ generated by
 $\phi_{\lm_i}$ and  its complex conjugate $\bar L_i \subset
 L^2(X)$.

    Then on the space $E$ the Hermitian form $H_{\bar L_i}$
 coincides with the form $H^{aut}_{\lm_i}$ corresponding to the
 automorphic trilinear form
  $l = l^{aut}_{\pi,\pi',\pi_i}: V \otimes V' \otimes V_{\lm_i} \to
  \bc$.

\end{claim}
   Indeed, if we identify the space $\bar L_i$ with $L_i^*$
   then the operator $pr_{\bar L_i} \circ r_\Dl: E \to \bar L_i$
    coincides
   with the operator $T^l$ corresponding to the automorphic
   trilinear form $l = l^{aut}_{\pi, \pi', \pi_i}$.

   This claim implies the first basic inequality

\begin{eqnarray}\label{spec-Q-Del-aut-1}
\sum_{\lm_i}|a_i|^2 H^{mod}_{\lm_i} \leq H_\Dl \ .
\end{eqnarray}

Indeed, by the uniqueness principle (\ref{coef-a-def}) we have:

\begin{eqnarray}\label{coef-a-Q}
H^{aut}_{\lm_i}=|a_i|^2 \cdot H^{mod}_{\lm_i}\ ,
\end{eqnarray}
where  $a_i=a_{\pi,\pi',\pi_i}$ are as in (\ref{coef-a-def}).

   Since all the spaces $\bar L_i$ are orthogonal we have
   $ \sum_i H^{aut}_{\lm_i} \leq H_\Dl$ which proves
   the first basic inequality.

\subsection{Second basic inequality }

   We would like to use the inequality (\ref{spec-Q-Del-aut-1})
 to bound the coefficients $a_i$. In order to do this we have
 to establish some bounds for the diagonal form $H_\Dl$.

The group $G\times G$ naturally acts on the space of Hermitian
forms on $C^\8(X\times X)$ - we denote this action by $\Pi$. We
extend this action to the action of the algebra
$\textsl{H}(G\times G) = C_c^\8(G\times G, \br)$ of smooth real
valued functions with compact support. Note that if $h \in
\textsl{H}(G\times G)$ is a nonnegative function then the operator
$\Pi(h)$ preserves the cone of positive forms.

We have then the second basic inequality

\begin{claim}{claimsecond}
Let $h \in \textsl{H}(G\times G)$ be a non-negative function. Then
there exists a constant $C$, depending on $h$, such that  we have
$\Pi(h) H_\Dl \leq C\cdot P_{X \times X}$, where $P_{X \times X}$
is the standard $L^2$ Hermitian form on the space $C^\8 (X \times
X)$.

\end{claim}

\begin{proof} Let $u \in C^\8(X \times X)$ and $f = |u|^2$.
Then $P_{X \times X}(u) = \la \mu ,f \ra $
 and $\Pi(h)H_\Dl(u) = \la \mu', f \ra$,
  where $\mu = \mu_{X \times X}$
  and $\mu' = \Pi(h)(\Dl_*(\mu_X))$
  are two measures on $X \times X$.

  Since the measure $\mu'$ is smooth it is bounded by $C \mu$.
  \end{proof}

 Note that the bound in the claim is essentially tight. Namely if the
function $h$ has large enough support then we also have a bound in
the opposite direction.

%%%%%%%%%%%%%%%%%%%%%%%%%%%%%%%%%%%%%%%%%%%%%%%%%%%%%%%%%%

\subsection{Positive functionals} \label{p-growth}
We can now prove that the coefficients $a_i$ have at most
polynomial growth in $|\lm_i|$.

We start with the inequality (\ref{spec-Q-Del-aut-1}) of {\it
non-negative forms}. We want to produce out of it an inequality
for coefficients $a_i$. There is a standard way to do this by
means of positive functionals on the space of Hermitian forms
$\CH(E)$.

 {\bf Definition.} A {\bf positive functional} on the space $\CH(E)$ is an
 additive map $\ro : \CH(E)^+ \to \br^+ \bigcup \8$.

  It is easy to see that the positive functional $\ro$ is automatically
  monotone and homogeneous (i.e. $\ro(H) \leq \ro(H')$ if $H \leq H'$ and
  $\ro(aH) = a\ro(H)$ for $a > 0$).

   \exa {example} Any vector $u \in E$ gives us an elementary
positive functional $\ro_u$ defined by $\rho_u(H) = H(u)$.

   Fix  a positive functional $\ro$ and consider the weight function
   $h(\lm)= \ro(H^{mod}_\lm)$. Then from the
   first basic inequality (\ref{spec-Q-Del-aut-1}) we can deduce the
    following inequality for a weighted sum of coefficients
    $|a_i|^2$:

$$ \sum_i h(\lm_i) |a_i|^2  \leq \ro (H_\Dl)\ .$$

\subsection{Test functional $\ro_T$}\label{testfunctional}
For any real $T$ we construct in  \ref{test} the positive "test"
functional $\ro_T$ on $\CH(E)$ with the properties described in
the proposition below. Let us fix automorphic representations $ V,
V' \subset C^\8(X), E = V \otimes V' \subset C^\8(X \times X)$ as
above.

\begin{prop}{rho-prop}
We can find a constant $C$ which depends only on $G$ and $\G$ and
a constant $T_0$ which depends on $V$ and $V'$ such that for any
$T\geq T_0$ there exists a positive functional $\ro_T$ on $\CH(E)$
satisfying
\begin{eqnarray}
&&\ro_T(H_\Dl)\leq CT^2\ , \label{ro-H_D}\\
&&h_T(\lm):=\ro_T(H_\lm^{mod})\geq 1  \ \text{ for any }
 |\lm|\leq 2T\ .\label{h(lm)}
\end{eqnarray}
\end{prop}

\subsection{Proof of Theorem \ref{THM}}

    Consider the inequality $ \sum_i |a_i|^2
\ro_T(H^{mod}_{\lm_i}) \leq \ro_T (H_\Dl)$.

The right hand side $\ro_T(H_\Dl)$ is bounded by $C T^2$. In the
left hand side we can leave only terms with $|\lm_i|\leq 2 T$.
Thus we arrive at inequality

\begin{eqnarray}\label{mean}
 \sum_{|\lm_i|\leq 2T}|a_i|^2\ \leq C T^2.
\end{eqnarray}

This gives the desired bound for $\sum |a_i|^2$.

According to Proposition \ref{propklambda}  there exists a
constant $b$ such that
 $b_iT^2\leq b |a_i|^2$
for $ T \leq |\lm_i| \leq 2T$. This shows that
$\sum_{T < |\lm_i| < 2T} b_i \leq A$
for some constant $A$, which finishes the proof of Theorem \ref{THM}.

\subsubsection{A conjecture} \label{a-conj}
One can show (see \cite{Re1}) that the mean-value result in
(\ref{mean}) is essentially sharp. One expects that for
$ T\leq |\lm_i|\leq 2T$
all terms in the sum (\ref{mean})  are at most of order
$T^\eps$ for any $\eps>0$.
Hence, we have established a sharp bound on the average
and a rather weak bound for each term.
This is  a typical situation
which one often encounters in the analytic theory of $L$-functions,
the so-called convexity bound.
The major problem hence is to find a method which would allow us
to obtain better bound
for a single term or for a short interval -- the
so-called subconvexity bounds.

We would like to make the following conjecture
concerning the size of coefficients $a_{\pi,\pi',\pi_i}$
which is equivalent to  Conjecture \ref{CONJ1}:

\begin{conj}{CONJ2}
 For fixed $\pi$, $\pi'$ and for any
$\eps>0$ there
exists $C_\eps>0$ independent of $\lm_i$ such that
$$|a_{\pi,\pi',\pi_{\lm_i}}|\leq C_\eps|\lm_i|^\eps\ ,$$
as $|\lm_i|\to\8$.
\end{conj}
%%%%%%%%%%%%%%%%%%%%%%%%%%%%%%%%%%%%%%%%%%%%%%%%%%%%%%%%%%%%%%%

\section{Construction of model trilinear functionals and of test functionals}\label{sect2}
\subsection{Model trilinear functionals} \label{modfunc}

For every $\lm \in \bc$ we denote by $(\pi_\lm,V_\lm)$  the smooth
class one representation of the generalized principle series of
the group $G=\PGLR$ described in \ref{irrrep}. We will use the
realization of $(\pi_\lm, V_\lm)$ in the space of smooth
homogeneous functions on $\br^2 \setminus 0$ of homogeneous degree
$\lm - 1$ .

For explicit computations it is often convenient to pass from
plane model to a circle model. Namely, the  restriction of
functions in $V_\lm$ to the unit circle $S^1 \subset \br^2$
defines an isomorphism of the space $V_\lm$ with the space
$C^\8(S^1)^{even}$ of even smooth functions on $S^1$ so we can
think about vectors in $V_\lm$ as functions on $S^1$.

In this section we describe the {\bf model} invariant trilinear
functional using the geometric models. Namely for given three
complex numbers $\lm_j$, $j=1, 2, 3$, we construct explicitly
nontrivial trilinear functional $\  l^{mod}: V_{\lm_1} \otimes
V_{\lm_2} \otimes V_{\lm_3} \to \bc$ by means of its kernel.

\subsubsection{Kernel of  $l^{mod}$ }\label{kernel}

 Let
$\om(\xi,\eta)=
\xi_1\eta_2-\xi_2\eta_1$ be $\SLR$-invariant of a pair of
 vectors $\xi,\ \eta\in\br^2$. We set

\begin{equation}\label{kern}
K_{\lm_1,\lm_2,\lm_3}(s_1,s_2,s_3)= |\om(s_2,s_3)|^{(\al-1)/2}\
 |\om(s_1,s_3)|^{(\beta-1)/2}|\om(s_1,s_2)|^{(\g-1)/2}
\end{equation}
for $s_1,s_2,s_3 \in\br^2\sm 0$, where $ \al =\lm_1-\lm_2-\lm_3,
  \beta= -\lm_1+\lm_2-\lm_3,\g = -\lm_1-\lm_2+\lm_3 $.

 The kernel function $K_{\lm_1,\lm_2,\lm_3}(s_1,s_2,s_3)$
satisfies two  main properties:
\begin{enumerate}
\item $K$ is invariant with respect to the diagonal action
of $\SLR$.

\item $K$ is homogeneous of degree $-1 -\lm_j$ in each variable $s_j$.
\end{enumerate}

 Hence
if $f_j$ are homogeneous functions of degree $-1+\lm_j$, then the
function
$$F(s_1,s_2,s_3)=
f_1(s_1)f_2(s_2)f_3(s_3)K_{\lm_1,\lm_2,\lm_3}(s_1,s_2,s_3)\ ,$$
is homogeneous of degree $-2$ in each variable $s_j\in\br^2\sm 0$.

\subsubsection{Functional $l^{mod}$ }\label{int-of-kern}
To define the model trilinear functional $l^{mod}$ we notice that
on the space $\CV$ of functions of homogeneous  degree $-2$ on
$\br^2\sm 0$ there
 exists a natural $\SLR$-invariant
functional $\fL:\CV \to\bc$ . It is given by the formula $\fL(f) =
\int_\Sigma f d\sigma$ where the integral is taken over any closed
curve $\Sigma\subset \br^2\sm 0$ which goes around $0$  and the
measure $d\s$  on $\Sigma$ is given by the area element inside of
$\Sigma$ divided by $ \pi$; this last normalization factor is
chosen    so that $\fL(Q^{-1}) = 1$ for the standard quadratic
form $Q$ on $\br^2$.

Applying $\fL$ separately to each variable $s_i\in \br^2\sm 0$ of
the function $F(s_1,s_2,s_3)$ above we obtain the $G$-invariant
functional
\begin{equation}\label{pairing}
l^{mod}_{\pi_1,\pi_2,\pi_3}(f_1\otimes f_2\otimes f_3):= \langle
\fL\otimes\fL\otimes\fL ,F \rangle\ .
\end{equation}

In the circle model this functional is expressed by the following
integral:
\begin{equation}\label{circleintegral}
l^{mod}_{\pi_1,\pi_2,\pi_3}(f_1\otimes f_2\otimes f_3)= (2
\pi)^{-3} \iiint f_1(x)f_2(y)f_3(z)K_{\lm_1,\lm_2,\lm_3} (x,y,z)
dx dy dz ,
\end{equation}

where $x, y , z$ are the standard angular parameters on the
circle.

\rem {remark} The integral defining the trilinear functional is
often divergent and the functional should be defined using
regularization of this integral. There are standard procedures how
to make such a regularization (see e.g. \cite{G1}).

   Fortunately in the case of unitary representations
    all integrals converge absolutely so we will not discuss the
regularization procedure.

%%%%%%%%%%%%%%%%%%%%%%%%%%%%%%%%%%%%%%%%%%%%%%%%%%%%%%%%%%%%

\subsection{Construction of  test functionals}\label{test}

 In this section we
will present a construction of a family of test functionals
$\ro_T$ on the space $\CH(E)$.

 Fix smooth irreducible pre-unitary representations
of class one $V = V_\tau, V = V_{\tau'}$ and denote by $E$ the
smooth completion of $V \otimes V'$ as in section \ref{pbound}. We
will do the computations only for representations of the principal
series; complementary series are treated similarly.

   For computations we will identify the spaces $V$ and $V'$ with
   $C^\8(S^1)^{even}$.

 Our aim is to prove the following

\begin{prop}{proptest}
There exist constants $T_0, C,\ c>0$ such that for any $T\geq T_0$
there exists a positive functional $\ro$ on $\CH(E)$ satisfying
\begin{eqnarray}
&&\ro(H_\Dl)\leq CT^2\ , \label{ro-H_D-app-B}\\
&&h_T(\lm):= \ro(H_\lm^{mod})\geq c \ \text{for any } |\lm|\leq
2T\ . \label{h(lm)-app-B}
\end{eqnarray}
\end{prop}

 The functional $\ro_T = c^{-1} \ro$ is the required test
 functional in \ref{testfunctional}.

\subsection{Proof of Proposition \ref{proptest}}\label{prooftest}

We will construct a functional $\ro$ as an integral of elementary
functionals. Namely we find  a positive function $h \in
\textsl{H}(G\times G) \subset C^\8(G \times G)$  and a vector $u
\in E$
  and define
$\ro(H)= \ro_u (\Pi(h)(H))$, where $\ro_u$ is the elementary
functional on the space $\CH(E)$ corresponding to the vector $u$.

\subsubsection{Construction of function $h$} We construct the function $h$
independent of parameter $T$.
 Let $D_1 \subset \SLR  \subset G $
be the subset of matrices $g$ with $||g|| \leq 2$. We consider the
subset $D = D_1 \times D_1 \subset G \times G$ and choose a
positive function $h \in \textsl{H}(G \times G) = C_c^\8(G \times
G)$ which is $\geq 1$ on the subset $D$ and is supported in some
neighborhood of $D$. We also assume that the function $h$ is
invariant under left and right translations by elements of the
maximal compact subgroup $K \times K$.

  \subsubsection{Construction of vector $u$}\label{consttest}

    Let us identify the space $E = V \otimes V'$ with the space
    of smooth functions $C^\8(S^1 \times S^1)^{even}$. Let $S$ be a disc
    in $S^1 \times S^1$ of radius $(100 T)\inv$. We construct $u$ as
    a smooth non-negative real valued function on $S^1 \times S^1$ supported in $S$
    such that

    (i) $\int u dx dy = 1$,

    (ii) $||u||_{L^2}^2 \leq 10^5 T^2$.

    We would like to show that the functional $\ro$ constructed
    in \ref{prooftest}
    satisfies conditions formulated in  Proposition \ref{proptest}.

    \subsubsection{ Geometric bound}\label{geometricbound} We have
    $$\ro(H_\Dl) = \ro_u (\Pi(h) (H_\Dl)) \leq C' \ro_u (P_E)
    = C' P_E(u) \leq C T^2$$
     see \ref{claimsecond}, \ref{consttest}.

    \subsubsection{ Spectral bound}\label{spectralbound} First we would like to
    give  another description of the Hermitian form $H^{mod}_\lm$.
    Consider the model trilinear
    functional $l = l^{mod}_{\pi,\pi',\pi_i}$ described in \ref{modfunc}
    and the corresponding operator $T^l: E \to \bar V_\lm$.

    We will identify the space $V_\lm$ with the space $C^\8(S^1)^{even}$.
    Fix some point $z \in S^1$. Then the $\delta$ functional $\dl_z$
    at this point
    defines a functional $f$ on $E$. We denote by $P_f$ the
    corresponding Hermitian form on $E$, $P_f(u) = |\langle f,u \rangle |^2$.

    Since the scalar product on the space $\bar V_\lm \simeq
    C^\8(S^1)^{even}$ is given by the standard integral we see that
    the standard Hermitian form $P_{\bar V_\lm}$ is an average
    over the compact group $K$ of the forms $\pi(k)(P_{\dl_z})$.
    This implies that
 $H^{mod}_\lm = \int_K (\Pi(k,k)(P_f)) dk$.

   Since we assumed the function $h \in \textsl{H}(G \times G)$ to
   be $K \times K$-invariant we see that

   $$\ro(H_\lm^{mod}) = \ro_u (\Pi(h)(P_f))$$

   Thus we see that in order to prove a lower bound for
   $\ro(H_\lm^{mod})$ it is enough to establish a lower bound
   for $\ro_u (\Pi(g) P_f) : = |\langle \Pi(g)f,u\rangle |^2$,
   for a subset of $g\in D$ of a measure bounded from below by a constant.

   Namely, the desired lower bound follows from the following

   \begin{lem} {D0} Let $T_0$ be large enough. Then there exists an open
    non-empty subset $D_0 \subset D$ such that for $T \geq T_0$ and for
   $g \in D_0$ we have
  $|\langle \Pi(g) f ,u \rangle | \geq 1/2$.
\end{lem}
   \begin{proof} Let $x,y$ be parameters on circles $S^1$ describing $V$
and $V'$. As follows from the definition (\ref{circleintegral})
the functional $f$ is given by the function
    $f =f(x,y)$  on $S^1 \times S^1$ described by
     $$f(x,y) =|\sin(y-z)|^{(\al-1)/2}
     |\sin(x-z)|^{(\beta-1)/2} |\sin(x-y)|^{(\g-1)/2}\ ,
     $$
where $\al,\ \beta,\ \g \in i\br$.

   Let $D_0 \subset D$ be the subset of  elements $g \in D$
   such that the restriction of $\Pi(g)(f)$  to the  subset $S
   \subset S^1 \times S^1$ has the absolute value $\leq 10$. (Note that
   the absolute value $|\Pi(g)(f)|$ is bounded from {\it below} for any $g\in D$
   by a constant depending only on $D$.)
   It is easy to see that for large $T$ the set $D_0$ is a non empty subset of
   $D\subset G\times G$ of a measure  bounded from below by a constant which
   is independent of $T$.

   On the other hand, for $g\in D_0$ we see that the gradient
   of the function $\Pi(g)(f)$ on the subset $S$ is bounded by $3T$.
   We note now that the diameter of $S$ is bounded by $(100T)\inv$
   and
   hence the lower bound on $|\langle \Pi(g)f,u \rangle |$ for $g \in D_0$
   follows from the following easy claim

   \begin {claim}{claim}Let $S$ be a set with a measure $\nu$ and $u, h$ be
   two measurable functions on  $S$. Let us assume that

   (i) $u$ is real valued positive function and $\int u d\nu = 1$.

   (ii) $sup|h(s)| \geq 1$  and the variation
   ${\text Var}(h) := sup |h(s) - h(s')|$ is bounded by $1/2$.

   Then $|\int h u d \nu| \geq 1/2$.

   \end{claim}

\end{proof}
%%%%%%%%%%%%%%%%%%%%%%%%%%%%%%%%%%%%%%%%%%%%%%%%%%%%%%%%%%%%%%%%%%%
\subsection{Construction of  test functionals via Sobolev
norms}\label{sobolev}

    In this section we outline another, slightly more conceptual,
construction of  test functionals. This construction uses the
notion of Sobolev norms on representation spaces (see \cite{BR2}).
\subsubsection{Sobolev norms}
Let $G$ be a Lie group and $(\pi, G, V)$ a smooth pre-unitary
representation. Then we can construct  a family of positive
definite Hermitian forms on the space $V$ as follows.

Fix a basis $\{ X_j |j = 1,...,r \}$ of the Lie algebra ${\mathfrak
g}$ of the group $G$. Then for any natural number $l$ and any $T
>0$ we define a Hermitian form $Q_{l,T}$ on V by
$$ Q_{l,T}(v) = \sum_\nu T^{2(l-|\nu|)}P(X^\nu(v))\ .$$
Here  the sum is over all multi indexes $\nu = (n_1,...,n_r)$ with
the norm $|\nu| := \sum n_j$ bounded by $l$ and $P = P_V$ is the
Hermitian form defining the unitary structure on $V$.

\subsubsection{Positive functionals defined by forms}

   Every positive definite Hermitian form $Q$ on $V$ defines a
positive functional $\ro_Q$ on $\CH(V)$ by $\ro_Q(H) = tr(H|Q)$.
   Here $tr(H|Q)$ denotes  the relative trace of forms $H$ and
   $Q$; by definition it is equal to the square of the
   Hilbert-Schmidt norm of the identity operator on $V$ considered
   as a morphism of pre-Hilbert spaces $(V,Q) \to (V,H)$.
This notion is discussed in detail in \cite{BR2}.

\subsubsection{Construction of  Sobolev test functionals}
Let us apply these constructions to the representation $(\Pi, G
\times G, E)$ discussed in \ref{pbound}.

 Fix $l$ and $T$, consider the Sobolev Hermitian form $Q = Q_{l,T}$
   on the space $E$ and define the positive functional $\ro$ on
$\CH(E)$ to be $\ro = \ro_{Q}$.

\begin{prop}{testsobolev} Suppose $l \geq 2$. Then

(i) $\ro(H_\Dl) \leq C T^{2-2l}$,

(ii) There exists $c > 0$  such that $\ro(H^{mod}_\lm) \geq c
T^{-2l}$ for $|\lm| \leq 2T$.

\end{prop}

   This gives another proof of Proposition \ref{testfunctional}.

   \subsubsection{Sketch of the proof of Proposition
   \ref{testsobolev}}

   (i) Since  the representation $\Pi$ is continuous with respect
   to the form $Q_{l,T}$ the second basic inequality \ref{claimsecond} implies
   that $\ro(H_\Dl) \leq C' \ro(P_E)$. The proof of the inequality
   $\ro(P_E) \leq C'' T^{2-2l}$ is the same as in \cite{BR2}, section 4.

   In order to prove (ii) it is enough to find a vector $u \in E$
   such that $Q_{l,T}(u) \leq T^{2l}$ and $|\langle f,u \rangle | \geq c$,
   where $f = f_z$ is the function described in \ref{spectralbound}.
   We can take a function $u \in C^\8(S^1 \times S^1)$  of the
   form $u = \phi f$ where $\phi$ is a smooth cut-off function which
   equals $0$ around singularities of the function $f$.

   We leave details to the reader.

%%%%%%%%%%%%%%%%%%%%%%%%%%%%%%%%%%%%%%%%%%%%%%%%%%%%%%%%%%%%%%%%%%%%%%%%%
\appendix
 \section { }    \label{A}
\subsection{Computation of $l^{mod}$ for $K$-fixed vectors}
\label{compkvect}

    In this appendix we prove the Proposition \ref{propklambda}
    which describes the assymptotic behavior of the function
    $k_\lm$.

One can prove this proposition applying the stationary phase
method directly to the integral (\ref{circleintegral}). To do this
we need to consider the complexification of the functions
$e_\lm(s_i)$ and the function $K_{\lm_1,\lm_2,\lm_3}(s_1,s_2,s_3)$
 in the variables $s_i$ and move contour of integration towards
the singularities of the complexified integral. This could be done
either in a classical language or using analytic continuation of
representations  in the spirit of \cite{BR1}.

\subsection{Computation of the integral}

   We prefer to prove this proposition in a different way. Namely we
   explicitly compute the value of the model functional
   on the unit vectors in terms of $\G$-
   functions and then prove the proposition
   by applying Stirling formulas for assymptotic behavior of
   $\G$-functions.

Let $\pi_{\lm_i}$, $i=1,2,3$ be three representations of the
generalized principal series and $e_{\lm_i}$ be the corresponding
$K$-fixed  unit vectors (they correspond to function $1$ in the
circle model).
   Set $A(\lm_1,\lm_2,\lm_3): =
   l^{mod}_{\pi_{\lm_1},\pi_{\lm_2},\pi_{\lm_3}} (e_{\lm_1}\otimes
e_{\lm_2}\otimes e_{\lm_3})$.

   In sections  \ref{B}, \ref{C}  we explicitly compute the function
   $A(\lm_1,\lm_2,\lm_3)$ (see the final expression in \ref{C}) .

\subsection{Gaussian}        \label {Gaussian}

   We would like to compute our integral by comparing it with
   Gaussian integrals which are much easier to manipulate with.

   Namely, suppose we are given a finite-dimensional
   Euclidean vector space $L$. In this case
   we introduce the Gaussian probability  measure $G$ on $L$ by
   $dG =  \pi^{- \dim L/2} \exp (-Q) dl$, where $Q$ is the quadratic
   form which defines the Euclidean structure on $L$ and $dl$ is the
   standard Euclidean measure on $L$.

      We are interested in the quantities
      $\langle f,G \rangle =\langle f,G \rangle_L:= \int f dG$ for
      various (usually homogeneous) functions $f$ on $L$. The main
      properties of the Gaussian which we use are the following:

      (i) {\bf Normalization}. $\langle 1 , G\rangle \ = 1$.

      (ii) {\bf Product formula}. Suppose that the Euclidean space
      $L$ is a product of Euclidean spaces $L_1$ and $L_2$. Then
      the Gaussian measure $G$ on $L$ is the product of Gaussian
      measures $G_1$ and $G_2$ on $L_1$ and $L_2$. In particular,
      if a function $f$ decomposes as a product of functions $f_1$
      and $f_2$ on $L_1$ and $L_2$ we have
      $\langle f,G\rangle = \langle f_1,G_1 \rangle \langle f_2,G_2 \rangle $.

   The following integrals are classical

\begin{prop}{classic}

    Let $L = \br^n$ be the standard Euclidean space.

      (i) Let  $r$ denote the radius
    function on $L$.
    Then $\la r^s,G \ra = \Gamma((s + n)/2)/\Gamma(n/2)$.

  (ii) Let $h$ be a linear functional on $L$. Then
  $\la |h|^s, G\ra = ||h||^s \Gamma((s + 1)/2)/\Gamma(1/2)$.

  (iii) Let $L$ be the space $M_{2,2}$ of  $\  2 \times 2$ matrices
  with the standard Euclidean structure. Then
  $\la |det|^s, G \ra = \Gamma((s+ 1)/2)\Gamma(s/2 + 1)
  /\Gamma(1/2)$.
\end{prop}

 \begin{proof}

  In (i) passing to spherical coordinates we get the integral
  $$2c \int r^{s+n-1} exp(- r^2)dr = c \int u^{(s+n)/2}exp(-u) du/u =
   c \Gamma((s + n)/2)\ .$$
    The normalization at $s = 0$
  defines the constant $c=1/\G(n/2)$.

  The proof of (ii) is reduced to the one variable case using
  product formula and then it follows from (i).

  In (iii) we can write $L$ as a product of two column spaces
  $L_1$ and $L_2$. Then we have

  $\la |det|^s,G \ra = \int |\omega(x,y)|^s dG_1(x) dG_2(y) =
  \int (\int|\omega(x,y)|^s dG_1(x))dG_2(y) =$

  $\Gamma((s + 1)/2)/\Gamma(1/2)\cdot \int|y|^s dG_2(y) =
  \Gamma((s+1)/2) \Gamma(s/2 +1)/ \Gamma(1/2)$

  since $\Gamma(1) = 1$.
\end{proof}

\subsection{Reduction 1}{\label{B}}
    The Proposition \ref{classic} allows us to write the integrals
    which we would like to compute as some Gaussian integrals.

\begin{cor}{corB}
  For any function $h \in V_{-\lm}$ we have
$$\la h,G\ra = \Gamma((1 -\lm)/2)\cdot \fL (h \cdot e_\lm)\ .$$
\end{cor}

   Indeed, after averaging $h$ with respect to the action of $SO(2)$
   we can assume that it is proportional to the function
   $e_{-\lm}$. Then the formula follows from \ref{classic}(i).

    Using this corollary we can rewrite the integral for the function
$A(\lm_1,\lm_2,\lm_3)$.

\begin{prop}{propAB} Consider the Euclidean space $L =  \br^2 \times \br^2
\times\br^2$ and define the function $B(\lm_1,\lm_2,\lm_3)$ by
Gaussian integral $B(\lm_1,\lm_2,\lm_3) := \la
K_{\lm_1,\lm_2,\lm_3}(s_1,s_2,s_3),G \ra $. Then
$$B(\lm_1,\lm_2,\lm_3) = A(\lm_1,\lm_2,\lm_3)\cdot {\G
((1-\lm_1)/2) \G((1-\lm_2)/2) \G((1-\lm_3)/2)}\ .$$
\end{prop}

\subsection{Reduction 2}{\label{C}}

  Let us rewrite the  integral defining the function $B$. First, we identify the
  Euclidean space $L$ in section \ref{B} with
  the space  $ M_{2,3}$ of $2 \times 3$ matrices. We consider the Euclidean space
  $W \approx \br^3$ and define the map $\nu: M_{2,3} \to W$ using
  $2 \times 2$ minors. Let us define the function $f$ on $W$ by
  formula $$f(w_1,w_2,w_3) =  |w_1|^{(\al
  -1)/2}|w_2|^{(\beta-1)/2}|w_3|^{(\g-1)/2}\ .$$
We can write $K_{\lm_1,\lm_2,\lm_3} = \nu^*(f)$ and hence
  $B(\lm_1,\lm_2,\lm_3) = \la \nu^*(f),G\ra $

  (here $ \al =\lm_1-\lm_2-\lm_3,
  \beta= -\lm_1+\lm_2-\lm_3, \g = -\lm_1-\lm_2+\lm_3$
   as in \ref{kernel}).

     Now we will use the following general lemma which we prove in
  section \ref{proofofcomparison}.

  \begin{lem}{comparison} Let $h$ be a homogeneous function on the space $W$
  of homogeneous degree $s$. Then
  $$\la \nu^*(h),G \ra_L = \la h, G\ra_W \cdot \G ( s/2 + 1  )\ .$$
\end{lem}

     From this lemma we see that the computation of the function
     $B(\lm_1,\lm_2,\lm_3)$ is reduced to the computation of the
     function $C(\al, \beta, \g):= \la f, G\ra $.

     Since the Gaussian $G$ on $W$ is a direct product of three
     one dimensional Gaussians and the function $f$ is a product
     of functions depending only on one coordinate we deduce that
     the integral $\la f,G\ra$ is a product of three one dimensional
     integrals which can be computed using \ref{classic}.

     Thus we obtain $C(\al, \beta, \g) = \G((\al+1)/4) \G((\beta+1)/4)
     \G((\g+1)/4)/\G(1/2)^3$.

     The final expression for the function $A(\lm_1,\lm_2,\lm_3) $ is

     $$A(\lm_1,\lm_2,\lm_3) = \frac { \G((\al+1)/4) \G((\beta+1)/4)
     \G((\g+1)/4)\G((\delta+1)/4)}
     {\G(1/2)^3\G((1-\lm_1)/2) \G((1-\lm_2)/2) \G((1-\lm_3)/2)} \ ,$$
where $\al =\lm_1-\lm_2-\lm_3, \beta= -\lm_1+\lm_2-\lm_3,
 \g = -\lm_1-\lm_2+\lm_3, \delta = -\lm_1-\lm_2-\lm_3$.

\subsection{Proof of Lemma \ref{comparison}} \label{proofofcomparison}

  Consider the natural actions of the group $SO(3)$ on the Euclidean spaces
  $M_{2,3}\approx W \times W$
  and $W$; these actions preserve Gaussian measures.

  The map $\nu: M_{2,3} \to W$ is $SO(3)$-equivariant; it is nothing else
  than
  the exterior product map $W \times W \to \bigwedge^2(W) = W^* =W$.
  Hence we can replace the
  function $h$ by its average with respect to the action of the group $SO(3)$,
  i.e. up to some constant by a function $h = r^s$.
  This shows that $\la \nu^*(h),G\ra = a(s) \la h,G\ra $, where
  $a(s)$ depends  on $s$ but not on $h$.

  In order to compute the function $a(s)$ we can consider the
  identity above for the function $h(w) = |w_3|^s$.
  According to Proposition \ref{classic}(ii) we have
   $\la h,G \ra = \G((s+1)/2)/ \G(1/2)$.

  On the other hand it is clear that the function $\nu^*(h)$
  depends only on four variables and hence the integral $\la \nu^*(h),G\ra $
  coincides with the integral $\la h',G \ra $ over the space $M_{2,2}$ of $2
  \times 2$ matrices , where $h'(m) = |\det(m)|^s$.

  From Proposition \ref{classic}(iii) we deduce that
  $a(s) = \G(s/2 + 1) $.

\subsection{Proof of Proposition \ref{propklambda}}
   According to Stirling formulas for any fixed $\sigma$ and large
   $t$ we have
   $\G(\sigma +i t) =\sqrt{2\pi} \exp (-\frac{\pi}{2} |t|)|t|^{\sigma -1/2}(1
   +O(|t|^{-1}))$.

   This and the explicit formula for the function
   $A(\lm_1,\lm_2,\lm_3)$ implies the proposition.

%%%%%%%%%%%%%%%%%%%%%%%%%%%%%%%%%%%%%%%%%%%%%%%%%%%%%%%%%%%%%%%%%

\end{document}